\documentclass[tbtags,reqno]{amsart}
  
\makeatletter
\newdimen\hoogte    \hoogte=12pt    
\newdimen\breedte   \breedte=14pt   
\newdimen\dikte     \dikte=0.5pt    
\def\beginYoung{
       \begingroup
       \def\vr{\vrule height0.8\hoogte width\dikte depth 0.2\hoogte}
       \def\fbox##1{\vbox{\offinterlineskip
                    \hrule height\dikte
                    \hbox to \breedte{\vr\hfill##1\hfill\vr}
                    \hrule height\dikte}}
       \vbox\bgroup \offinterlineskip \tabskip=-\dikte \lineskip=-\dikte
            \halign\bgroup &\fbox{##\unskip}\unskip  \crcr }
\def\End@Young{\egroup\egroup\endgroup}
\newenvironment{Young}{\beginYoung}{\End@Young}
\makeatother

\DeclareMathOperator{\Res}{Res}

\newtheorem{theorem}{Theorem}
\newtheorem{lemma}[theorem]{Lemma}

\newtheorem{conjecture}[theorem]{Conjecture}

\theoremstyle{remark}
\newtheorem*{acknow}{Acknowledgments}

\begin{document}

 \title[Creation operators for the Macdonald and Jack polynomials
 ]{Creation operators for the Macdonald and Jack polynomials} 
\author{Luc~Lapointe and Luc~Vinet}
\address
{Centre de Recherches Math{\'e}matiques, \\Universit{\'e} de
Montr{\'e}al,  C.P.~6128,  succursale~Centre-ville, \\ Montr{\'e}al,
Qu{\'e}bec, Canada, H3C 3J7}
 \date{}

 \begin{abstract}
   Formulas of Rodrigues-type for the Macdonald polynomials are
   presented.  They involve creation operators, certain properties of
   which are proved and other conjectured.  The limiting case of the
   Jack polynomials is discussed.
\end{abstract}
\maketitle

\section{Introduction}
A formula that gives the Jack polynomials through the application of a
string of creation operators on the Jack polynomials of lowest degree
has been introduced in \cite{2,3}. Since the Jack polynomials are a
specialization of the Macdonald polynomials which involve two
parameters $q$,$t$, it was natural to expect that a Rodrigues formula
should also exist for Macdonald polynomials.  We derive such a formula
in this paper.  It should be pointed out however, that the
construction presented here is much different than the one to be found
in \cite{2,3}; in fact, the latter is not simply the limit of the
former.  In this connection, we also give in the form of conjecture,
the expression of the creation operators that should be the analogs
the operators constructed in \cite{2,3}.  If this conjecture is true,
it implies in particular, that the expansion coefficients of the
Macdonald polynomials in the monomial basis are polynomials in $q$,$t$
with integer coefficients. The action of these operators on arbitrary
Macdonald polynomials is quite elegant and proves useful to conjecture
that families of $N$ commuting operators can be constructed out of
them.  The limiting case of the Jack polynomials is presented at the
end and the connection is made with the operators given in \cite{2,3}.
The conjectures are seen to be valid in this case also.

\section{Definitions \cite{1}}
Symmetric polynomials are labelled by partition $\lambda$ of their
degree $n$, that is sequences $\lambda = (\lambda_1, \lambda_2, \dots
)$ of non-negative integers in decreasing order $\lambda_1 \ge
\lambda_2 \ge \dots$ such that $|\lambda| = \lambda_1 + \lambda_2 +
\dots = n$.  The number of non-zero parts in $\lambda$ is denoted
$\ell(\lambda)$. Let $\lambda$ and $\mu $ be two partitions of $n$.  In
the dominance ordering, $\lambda \ge \mu $ if $\lambda_1 + \lambda_2 +
\dots + \lambda_i \ge \mu _1 + \mu _2 + \dots + \mu _i$ for all $i$.  We can
associate a diagram to each partition $\lambda$.  The diagram is made
out of $\ell(\lambda)$ rows, labelled by the integer $i$, with
$\lambda_i$ squares in each one of them.  The squares are identified
by the coordinates $(i,j)$ $\in \mathbb Z^2$ with $i$, the row index,
increasing as one goes downwards and j, the column index, increasing
as one goes from left to right.  For example the diagram of
$(5,4,4,1)$ is

$$
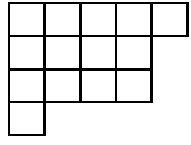
\begin{Young}
  & & & & \cr & & & \cr & & & \cr \cr
\end{Young}
$$

For each square $s=(i,j)$ in the diagram of a partition $\lambda$, let
$\ell'(s), \ell(s), a(s)$ and $a'(s)$ be respectively the number of
squares in the diagram of $\lambda$ to the north, south, east and west
of the square $s$.  By $\lambda \supset \mu $ it is meant that the
diagram ${\lambda}$ contains the diagram $\mu $, i.e. that $\lambda_i \ge
\mu _i$ for all $i \ge 1$.  The set-theoretic difference $\theta =
\lambda - \mu $ is called a skew diagram. For example if $\lambda =
(5,4,4,1)$ and $\mu =(4,3,2)$, $\lambda - \mu $ is made of the dotted
squares in the picture

$$
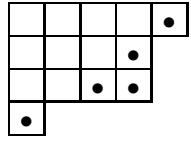
\begin{Young}
  & & & & $\bullet$ \cr & & & $\bullet$ \cr & & $\bullet$ & $\bullet$
  \cr $\bullet$ \cr
\end{Young}
$$

A skew diagram $\theta$ is a vertical $m$-strip if $|\theta|=m$ and
$\theta_i \le 1$ for each $i \ge 1$.  In other words, a vertical strip
has at most one square in each column.

Let $\Lambda_N$ denote the ring of symmetric functions in the
variables $x_1, x_2, \dots , x_N$.  Three standard bases for the space
of symmetric functions are:
\begin{enumerate}
\item[(i)] the power sum symmetric functions $p_{\lambda}$ which in
  terms of the power sums
\begin{equation}
p_i = \sum_k x_k^i,
\end{equation}
are given by
\begin{equation}
p_{\lambda} = p_{\lambda_1} p_{\lambda_2} \dots,
\end{equation}
\item[(ii)] the monomial symmetric functions $m_{\lambda}$ which are
\begin{equation}
m_{\lambda} = \sum_{\text{distinct permutations}} x_1^{\lambda_1} x_2^{\lambda_2} \dots
\end{equation}
\item[(iii)] the elementary symmetric functions $e_{\lambda}$ which in
  terms of the i-th elementary function
\begin{equation}
e_i = \sum_{j_1<j_2<\dots<j_i} x_{j_1} x_{j_2} \dots x_{j_i}= m_{(1^i)},
\end{equation}
are given by
\begin{equation}
e_{\lambda} = e_{\lambda_1} e_{\lambda_2} \dots .
\end{equation}
\end{enumerate}

The Macdonald polynomials can now be presented as follows.  To the
partition $\lambda$ with $m_i(\lambda)$ parts equal to $i$, we
associate the number
\begin{equation} \label{1}
z_\lambda
        = 1^{m_1} m_1 !  \, 2^{m_2} m_2! \dotsm
\end{equation}
Let $q$ and $t$ be parameters and $\mathbb Q(q,t)$ the field of all
rational functions of $q$ and $t$ with rational coefficients and
define a scalar product $\langle \ , \ \rangle_{q,t}$ on $\Lambda_N
\otimes \mathbb Q(q,t)$ by
\begin{equation}
\langle p_\lambda, p_\mu \rangle_{q,t}
        =\delta_{\lambda \mu } z_\lambda \prod_{i=1}^{\ell(\lambda)} \frac{1-q^{\lambda_i}}{1-t^{\lambda_i}},
\end{equation}
where $\ell(\lambda)$ is the number of parts of $\lambda$.  The
Macdonald polynomials $J_\lambda (x; q,t) \in \Lambda_N \otimes
\mathbb Q(q,t)$ are uniquely specified by
\begin{align}
  \mathrm{(i)} \ &  \langle J_\lambda, J_\mu \rangle_{q,t} = 0, \qquad \text{if } \lambda \ne \mu , \\
  \mathrm{(ii)} \ &  J_\lambda = \sum_{\mu \le \lambda} v_{\lambda\mu }(q,t) m_\mu , \\
  \mathrm{(iii)} \ & v_{\lambda\lambda}(q,t)= c_{\lambda}(q,t),
\end{align}
where
\begin{equation}
c_{\lambda}(q,t) = \prod_{s \in \lambda} (1-q^{a(s)} t^{\ell (s)+1}).
\end{equation}

For each $1 \le i \le N$ the shift operator $T_{q,x_i}$ is defined by
\begin{equation}
T_{q,x_i} f(x_1,\dots,x_N) = f(x_1,\dots,q x_i,\dots,x_N)
\end{equation}
for any polynomial $f \in \Lambda_N \otimes \mathbb Q(q,t)$.

For $r = 1,\dots,N$, let $M_N^r$ denote the Macdonald operator
\begin{equation}
 \quad M_N^r=\sum_I A_I (x;t) \prod_{i\in I} T_{q,x_i}
\end{equation}
summed over all $r$-element subsets $I$ of $\{1,\dots,N\}$, where
\begin{equation}
A_I (x;t) = t^{r(r-1)/2} \prod_{\begin{subarray}{c} i \in I \\ j \not \in I 
\end{subarray}} \frac{t x_i - x_j}{x_i - x_j},
\end{equation}
and take $M_N^0 \equiv 1$.  These operators commute with each other,
$[M_N^r,M_N^l]=0$ and are diagonal on the Macdonald polynomials basis.
From the Macdonald operators, one constructs
\begin{equation}
M_N(X;q,t) = \sum_{r=0}^N M_N^r X^r,
\end{equation}
with $X$ an arbitrary parameter.  The operator $M_N$ will play a
crucial role in the following.  Its action on $J_{\lambda}(x;q,t)$
with $\ell(\lambda) \le N$ is given, remarkably, by
\begin{equation}
M_N(X;q,t) J_{\lambda}(x;q,t) = a_{\lambda}(X;q,t) J_{\lambda}(x;q,t),
\end{equation}
where
\begin{equation}
a_{\lambda}(X;q,t) = \prod_{i=1}^N (1+X q^{\lambda_i} t^{N-i}).
\end{equation}
From (15) we see that the eigenvalue of $M_N^r$ on
$J_{\lambda}(x;q,t)$ is the coefficient of $X^r$ in the polynomial
(17).

Let us introduce the monic Macdonald polynomials $ P_{\lambda} =
1/c_{\lambda}(q,t) J_{\lambda}$.  The formulas giving the action of
the elementary symmetric function $e_k$ on the $P_{\lambda}$'s are
known as the Pieri formulas.  Explicitely, they are
\begin{equation}
e_k P_{\lambda} = \sum_\mu  \Psi_{\mu /\lambda} P_\mu 
\end{equation}
summed over partition $\mu \supset \lambda$ (of length $\le N$) such that
$\mu -\lambda$ is a vertical $k$-strip, with
\begin{equation} 
\Psi_{\mu /\lambda}= \prod_{\begin{subarray}{c} s \in C_{\mu /\lambda} \\ s 
\not \in R_{\mu /\lambda} \end{subarray}} \frac{b_\mu (s)}{b_{\lambda}(s)}
\end{equation}
and
\begin{equation}
b_{\lambda}(s) = {\cases \frac{1-q^{a(s)} t^{\ell(s)+1}}{1-q^{a(s)+1} t^{\ell(s)}} \qquad & {\text{if}}~s \in \lambda \\
1    & {\text{otherwise}} \endcases}
\end{equation}
where $C_{\mu /\lambda}$ (resp. $R_{\mu /\lambda}$) denote the union of the
columns (resp. rows) that intersect $\mu - \lambda$.  For example, with
$\mu =(4,2,2)$ and $\lambda = (3,2,1)$ we have

$$
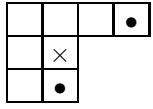
\begin{Young}
  & & & $\bullet$ \cr & $\times$ \cr & $\bullet$ \cr
\end{Young}
$$

\noindent so the only $s$ in   $C_{\mu /\lambda}$ but not in $R_{\mu /\lambda}$ is in  position (2,2).  With $\lambda+1$ representing the partition $(\lambda_1+1,\dots,\lambda_N+1)$, since $P_{\lambda+1}= e_N P_{\lambda}$, we have
\begin{equation}
 \Psi_{\mu +1/\lambda+1}=\Psi_{\mu /\lambda}.
\end{equation}

When restricted to the $N$-dimensional torus
$T=\bigl\{\{x_1,\dots,x_N\} \in \mathbb C^N; |x_i|=1, 1 \le i \le N
\bigr\}$, for $|q| <1$ and $|t|<1$, the polynomials $P_{\lambda}$ are
also orthogonal with respect to the scalar product $(~,~)$ defined by
\begin{equation}
(f,g)= \frac{1}{N} \int_T f(x) \overline {g(x)} \Delta(x;q,t) dx
\end{equation}
with $dx$ the normalized Haar measure on $T$ and
\begin{equation}
\Delta(x;q,t) = \prod_{i\neq j} (x_i x_j^{-1};q)_{\infty}/ (t x_i x_j^{-1};q)_{\infty},
\end{equation}
where
\begin{equation}
(a;q)_{\infty} = \prod_{i=0}^{\infty} (1-aq^i).
\end{equation}
In fact,
\begin{equation}
(P_{\lambda},P_\mu )=\delta_{\lambda \mu } c_N \prod_{s \in \lambda} b_{\lambda}(s)^{-1} \prod_{s \in \lambda} \frac{1-q^{a'(s)} t^{N-\ell'(s)}}{1-q^{a'(s)+1} t^{N-\ell'(s)-1}},
\end{equation}
where $c_N = (1,1)$.

\section{The Rodrigues formula}
\begin{theorem} The Macdonald polynomials $J_\lambda(x;q,t)$ associated to the partitions $\lambda = (\lambda_1 , \dots , \lambda_N)$ are given by
\begin{equation} 
J_\lambda(x;q,t)
        = (B_N^+)^{\lambda_N} (B_{N-1}^+)^{\lambda_{N-1} -\lambda_N} \dots (B_1^+)^{\lambda_1 - \lambda_2} \cdot 1, 
\end{equation}
with
\begin{equation}
B_k^+ = \frac{1}{(q^{-1};t^{-1})_{N-k}} M_N(-t^{k+1-N}q^{-1};q,t) e_k,\quad k=1,\dots,N,
\end{equation}
where for $n$ positive integer, $(a;q)_n= (1-a)(1-qa) \dots
(1-q^{n-1}a)$ and $(a;q)_0 \equiv 1$.
\end{theorem}
{\it{Proof.}} From the Pieri formula, we have the following lemma
\begin{lemma} The action of $e_k$ on $P_{\lambda}$ with $\lambda$ a partition  with  $\ell(\lambda) \leq k$  is given by 
\begin{equation}
e_k P_{\lambda} =  P_{\lambda+(1^k)} +
\sum_{\mu \neq \lambda+(1^k)} \Psi_{\mu /\lambda} P_\mu ,
\end{equation} 
where all the $\mu $'s in the sum are such that $\mu _{k+1}=1$.
\end{lemma}
This is seen from the fact that $\mu $ must be a partition which contains
$\lambda$ and $\mu -\lambda $ a vertical $k$-strip.  Hence the only way
to construct a $\mu $ with $\mu _{k+1} \neq 1$ is to add a 1 in each of the
first $k$ entries of $\lambda$.

From Lemma~2 and (16), we thus have
\begin{equation}
 M_N(-t^{k+1-N}q^{-1};q,t) e_k  P_{\lambda} = \prod_{i=1}^k (1-t^{k+1-i}q^{\lambda_i}) (q^{-1};t^{-1})_{N-k} P_{\lambda+(1^k)},
\end{equation}
since the eigenvalues of $ M_N(-t^{k+1-N}q^{-1};q,t)$ on the $P_\mu $'s
in (28) are
\begin{equation}
a_\mu (-t^{k+1-N}q^{-1};q,t) = \prod_{i=1}^N (1-t^{k+1-i}q^{\mu _i-1}),
\end{equation}
and vanish if $\mu _{k+1}=1$.
\begin{lemma} If $\lambda$ is a partition with  $\ell(\lambda) \leq k$ ,
\begin{equation}
\frac{c_{\lambda+(1^k)}}{c_{\lambda}}=\prod_{i=1}^k (1-t^{k+1-i}q^{\lambda_i}).
\end{equation}
\end{lemma}
When going from $\lambda$ to $\lambda+(1^k)$, what we actually do is
add a column at the west of the diagram.

$$
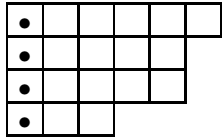
\begin{Young}
  $\bullet$ & & & & & \cr $\bullet$ && & & \cr $\bullet$ && & &\cr
  $\bullet$ &&\cr
\end{Young}
$$

\noindent From the fact that $c_{\lambda}$ only involves $a(s)$ and $\ell(s)$ which do not  depend on the number of square at their west, the contribution in $c_{\lambda+(1^k)}$ of the squares that have been shifted to the east is exactly $c_{\lambda}$.  Hence we only have to take the product of the contributions of the first column of $\lambda+(1^k)$, which is $\prod_{i=1}^k (1-t^{k+1-i}q^{\lambda_i})$.  

Using Lemma~3 and passing from $P_{\lambda}$ to $J_{\lambda}$, we have
\begin{equation}
\frac{1}{(q^{-1};t^{-1})_{N-k}} M_N(-t^{k+1-N}q^{-1};q,t) e_k  J_{\lambda}= J_{\lambda+(1^k)}
\end{equation}
which gives when $\ell(\lambda) \leq k$,$B_k^+ J_{\lambda}=
J_{\lambda+(1^k)}$ and hence, proves Theorem~1.

Note that the adjoint of $B_k^+$ with respect to the scalar product
$(~,~)$ is easily found since $M_N(X;q,t)$ is self-adjoint and
$e_k^{\dagger}=e_N^{-1}e_{N-k}$.  Indeed, one has
\begin{equation}
B_k^-=(B_k^+)^{\dagger}= \frac{1}{(q^{-1};t^{-1})_{N-k}} e_N^{-1}
e_{N-k} M_N(-t^{k+1-N}q^{-1};q,t) . 
\end{equation}
Using these operators and the fact that $(P_\mu ,P_{\lambda})=0$ if
$P_\mu \neq P_{\lambda}$, allows for a straightforward computation of
$(P_{\lambda},P_{\lambda})$.  One first observes that if
$\ell(\lambda) \leq k$,
\begin{equation}
\begin{split}
  (P_{\lambda+(1^k)},P_{\lambda+(1^k)})&=\frac{c_{\lambda}}{c_{\lambda+(1^k)}}
  (B_k^+ P_{\lambda},P_{\lambda+(1^k)})\\ 
  &=\frac{c_{\lambda}}{c_{\lambda+(1^k)}} ( P_{\lambda},B_k^-
  P_{\lambda+(1^k)}).
\end{split}
\end{equation}
From (31) and the eigenvalue of $M_N(-t^{k+1-N} q^{-1};q,t)$ on
$P_{\lambda+(1^k)}$ given in the R.H.S. of (29) , one sees that the
constant cancels out and that
\begin{equation}
\begin{split}
  (P_{\lambda+(1^k)},P_{\lambda+(1^k)})&=
  (P_{\lambda},e_N^{-1} e_{N-k} P_{\lambda+(1^k)})\\
  &=(P_{\lambda+1}, e_{N-k} P_{\lambda+(1^k)}).
\end{split}
\end{equation}
Finally, using the orthogonality of the $P_{\lambda}$'s and the Pieri
formula, one finds the following formula
\begin{equation}
(P_{\lambda+(1^k)},P_{\lambda+(1^k)})= \Psi_{\lambda+1/\lambda+(1^k)}
(P_{\lambda},P_{\lambda}), 
\end{equation}
which gives the norm of $P_{\lambda}$ through iteration.

\section{Conjectures}

The Rodrigues formula given in Theorem~1 does not imply that the
$v_{\lambda \mu }(q,t)$'s of (9) are polynomials in $q$ and $t$ with
integer coefficients.  However, it proves useful to obtain results
which once proved, would have this implication.  Such formulas which
represent generalizations for the Macdonald polynomials of relations
that we have proved for the Jack polynomials \cite{2,3}, are given
below in the form of conjectures.  Their limits as $ q=t^{\alpha}$ and
$t \to 1$ will be discussed in the section on the Jack polynomials.
\begin{conjecture} The Macdonald polynomials $J_\lambda(x;q,t)$  are given by
\begin{equation} 
J_\lambda(x;q,t)
        = (\tilde B_N^+)^{\lambda_N} (\tilde B_{N-1}^+)^{\lambda_{N-1}
          -\lambda_N} \dots (\tilde B_1^+)^{\lambda_1 - \lambda_2} \cdot 1,        
\end{equation}
with
\begin{equation}
\tilde B_k^+ = \sum_I \tilde A_I(x;t) x_I M_I(-t;q,t)
\end{equation}
summed over all $k$-element subsets $I$ of $\{1,\dots,N\}$, where
\begin{equation}
x_I = \prod_{i \in I} x_i,
\end{equation}
\begin{equation}
\tilde A_I(x;t) = t^{-(N-k)k} \prod_{\begin{subarray}{c} i \in I \\ j \notin I
\end{subarray}} \frac{t x_i -x_j}{x_i - x_j}
\end{equation}
and
\begin{equation}
M_I(X;q,t)=M_k (X;q,t) 
\end{equation}
in the $k$ variables $x_i \in I$.
\end{conjecture}
  
\begin{conjecture} The Macdonald polynomials $J_\lambda(x;q,t)$  are given by
\begin{equation} 
J_\lambda(x;q,t)
        = (\bar B_N^+)^{\lambda_N} (\bar B_{N-1}^+)^{\lambda_{N-1}
          -\lambda_N} \dots (\bar B_1^+)^{\lambda_1 - \lambda_2} \cdot 1,  
\end{equation}
with
\begin{equation}
\bar B_k^+ = \sum_I \bar A_I(x;t) x_I M_I(-t;q,t)
\end{equation}
summed over all $k$-element subsets $I$ of $\{1,\dots,N\}$, where
\begin{equation}
\bar A_I(x;t) = \biggl\{ \prod_{\begin{subarray}{c} i \in I \\ j \not \in I
\end{subarray}}  \frac{ x_i -t x_j}{x_i - x_j} \biggr\} \biggl\{
\prod_{  j \not \in I }  T_{q,x_j} \biggr\}. 
\end{equation}
\end{conjecture} 
 
A manifest corollary of Conjecture~5 would be that the $v_{\lambda
  \mu }(q,t)$'s are polynomials in $q,t$ with integer coefficients.  Let
us stress however that the operators $\tilde B_k^+$ of Conjecture~4
appear to be the natural generalizations of the creation operators
introduced in \cite{3} in the case of the Jack polynomials. Indeed, as
will be confirmed in section~6, the operators (38) and (85) share
important properties.

As an indication that Conjecture~4 must be true, we prove that, for
partition $\lambda$ with $\ell(\lambda) \leq k$,
\begin{equation}
\tilde B_k^+ J_{\lambda} =  J_{\lambda +(1^k)},
\end{equation}
in the cases where the number of variables $N=k$ or $k+1$.
\begin{lemma}
\begin{equation}
M_I(X;q,t) e_N = e_N M_I(X q ;q,t)
\end{equation}
for all subsets $I$ of $\{1,\dots,N\}$.
\end{lemma}
This result follows from the fact that $M_I^r e_N = q^r e_N M_I^r$,
which is easily derived from (13).

Lemma~6 immediately ensures that (45) is true for $N=k$, since it
implies that $\tilde B_k^+ = B_k^+ $ in this case.  The following
lemma will be needed in the proof of the special case $N=k+1$.
\begin{lemma}
  When the number of variables is $k+1$,
\begin{equation}
M_{k+1}(-1;q,t) J_{(\lambda_1,\dots,\lambda_k,0)}=0
\end{equation}
\end{lemma}
Since $a_{\lambda}(-1;q,t)= \prod_{i=1}^{k+1} (1-q^{\lambda_i}
t^{k+1-i})$, setting $\lambda_{k+1}=0$ implies that $a_{\lambda}=0$.
From (16), this lemma is then seen to hold.

Given Lemma~7, showing that
\begin{equation}
\tilde B_k^+ = B_k^+ +G M_{k+1}(-1;q,t),
\end{equation}
when $N=k+1$, with $G$ a certain expression, will prove (45).  The
following identity obtained by Garsia and Tesler \cite{4} will be used
to prove (48):
\begin{equation}
\sum_I \prod_{\begin{subarray}{c} i \in I \\ j \not \in I
\end{subarray}}  \frac{ x_i -t x_j}{x_i - x_j} x_I =\sum_I x_I= e_k,
\end{equation}
with the sum over all $k$-element subsets $I$ of $\{1,\dots,N\}$.

In the special case $N=k+1$,
\begin{equation}
B_k^+= (1-q^{-1})^{-1} \sum_{r=0}^{k+1}(-q^{-1})^r M_{k+1}^r e_k.
\end{equation}
Taking the part of this equation involving $M_{k+1}^r x_1\dots x_k$,
we have
\begin{equation}
\begin{split}
  M_{k+1}^r & x_1\dots x_k = \\
  & q^{r-1} x_1\dots x_k M_{k+1}^r + (q^r-q^{r-1})x_1\dots x_k
  \sum_{\begin{subarray}{c} I \subset \{1,\dots,k\} \\ |I|=r
\end{subarray}} A_I(x;t) \prod_{i \in I} T_{q,x_i}
\end{split}
\end{equation}
remembering that in $A_I(x;t)$ the total set of variables is
$\{1,\dots,k+1\}$.  The second term in the R.H.S. of (51) is
\begin{equation}
\begin{split}
  x_1\dots x_k \sum_{\begin{subarray}{c} I \subset \{1,\dots,k\} \\ 
      |I|=r
\end{subarray}} & A_I(x;t) \prod_{i \in I} T_{q,x_i}=\\
& \sum_{\begin{subarray}{c} I \subset \{1,\dots,k\};|I|=r \\ I^c =
    \{1,\dots,k\} \backslash I
\end{subarray}} x_{I^c} t^{r(r-1)/2}  \prod_{\begin{subarray}{c} i \in
  I \\ j  \not \in I 
\end{subarray}} \frac{ t x_i - x_{j}}{x_i - x_{j}} x_I \prod_{i \in I}
T_{q,x_i}. 
\end{split}
\end{equation}
With $I$ fixed, knowing that to compute $M_{k+1}^r e_k$ we will have
to sum over all $I'\subset \{1,\dots,k+1\} \backslash I$ with
$|I'|=k-r$, of which $I^c$ is a special case, and using the following
special case of the Garsia-Tesler formula:
\begin{equation}
 \sum_{\begin{subarray}{c} I' \subset \{1,\dots,k+1\}\backslash I\\
     |I'|=k-r \end{subarray}} x_{I'} =  \sum_{\begin{subarray}{c} I'
     \subset \{1,\dots,k+1\}\backslash I\\ |I'|=k-r \end{subarray}}
 \prod_{\begin{subarray}{c} i \in I' ;j \not \in I' \\ j \in
     \{1,\dots,k+1\}\backslash I \end{subarray}} \frac{  x_i - t^{-1}
   x_{j}}{x_i - x_{j}} x_{I'}, 
\end{equation}
we see that we can replace $x_{I^c}$ by $t^{r-k} \prod_{i \in I^c}
\frac{ t x_i - x_{k+1}}{x_i - x_{k+1}} x_{I^c}$ in (52).  Equation
(52) thus becomes
\begin{equation}
\begin{split}
  x_1 & \dots x_k \sum_{\begin{subarray}{c} I \subset \{1,\dots,k\} \\ 
      |I|=r
\end{subarray}} A_I(x;t) \prod_{i \in I} T_{q,x_i}=\\
& t^{r-k} \prod_{i \in \{1,\dots,k\}} \frac{ t x_i - x_{k+1}}{x_i -
  x_{k+1}} x_1 \dots x_k \sum_{\begin{subarray}{c} I \subset
    \{1,\dots,k\};|I|=r \\ I^c = \{1,\dots,k\} \backslash I
\end{subarray}} t^{r(r-1)/2} \prod_{\begin{subarray}{c} i \in I \\ j  \in I^c
\end{subarray}}  \frac{ t x_i - x_{j}}{x_i - x_{j}}  \prod_{i \in I} T_{q,x_i}.
\end{split}
\end{equation}
Hence, from (51) and (54), the part of $B_k^+$ in (50) involving $ x_1
\dots x_k $ is
\begin{equation}
q^{-1} (1-q^{-1})^{-1}  x_1 \dots x_k M_{k+1}(-1;q,t)+ x_1 \dots x_k
\tilde A_{\{1,\dots,k\}}(x;t) M_{\{1,\dots,k\}}(-t;q,t), 
\end{equation}
which by symmetry implies (48).

We see that $\tilde B_k^+$ and $B_k^+$ coincide only on the set of
Macdonald polynomials with $\ell(\lambda) \le k$.  The action of
$\tilde B_k^+$ on an arbitrary Macdonald polynomial will be given
below in the form of a conjecture.  It is much simpler than that of
$B_k^+$.  The operators $\tilde B_k^+$ further have a number of
remarkable properties that the operators $B_k^+$ do not share.

\section{ Properties of the operators $\tilde B_k^+$}

We extend the definition of a partition to allow real entries:
\begin{equation}
P_{(\beta_1,\dots,\beta_N)} = e_N^{\beta_N}
P_{(\beta_1-\beta_N,\dots,\beta_{N-1}-\beta_N,0)}= e_N^{\beta_N}
P_{\beta-\beta_N} 
\end{equation}
$\forall \beta_N \in \mathbb R$ and $\beta_{i}-\beta_{i+1}$ an integer
$\ge 0$, $i=1,\dots,N-1$.  Take the operator $F(\kappa)$ that acts as
follows on $P_{\beta}$:
\begin{equation}
F(\kappa) P_{(\beta_1,\dots,\beta_N)}=\prod_{i=1}^N (t^{\kappa+1-i}
q^{\beta_i-1};q^{-1})_{\infty} P_{(\beta_1,\dots,\beta_N)}, 
\end{equation}
From this definition and upon defining $q=t^{\alpha}$, it follows that
\begin{equation}
F(\kappa) e_N^{\rho}=e_N^{\rho} F(\kappa+\alpha \rho).
\end{equation}
We now form the operators
\begin{equation}
F_{m,\kappa} = F(\kappa) e_m F(\kappa)^{-1},
\end{equation}
to see that these have on $P_{\beta}$ actions that only involve a
finite number of products.  These are given by
\begin{equation}
F_{m,\kappa} P_{\beta} = \sum_{\delta} \Psi_{\delta/\beta}
F_{\delta/\beta}(\kappa) P_{\delta}
\end{equation}
with $\delta-\beta$ $m$-vertical strips,
\begin{equation}
F_{\delta/\beta}(\kappa)= \prod_{\begin{subarray}{c} s \in
    \delta-\beta_N \\ s \not \in \beta-\beta_N 
\end{subarray}} F_{\delta-\beta_N}(s;\kappa+\alpha \beta_N),
\end{equation}
and where for partitions $\mu $ made of nonnegative integers,
\begin{equation}
F_\mu (s;\kappa)=(1-t^{\kappa-\ell'(s)} q^{a'(s)}),\quad \forall s \in \mu .
\end{equation}
Note that by the same argument as in (21), $\Psi_{\delta/\beta} =
\Psi_{\delta-\beta_N/\beta-\beta_N}$.  Hence, the action of
$F_{m,\kappa}$ on the Macdonald polynomials is very similar to that of
$e_m$ on these functions: the action of $F_{m,\kappa}$ differs from
that of $e_m$ by the presence of $m$ additional factors in front of
the coefficients $\Psi_{\delta/\beta}$.  As an example, taking $N=4$
and $m=2$,
\begin{equation}
\begin{split}
  F_{2,\kappa} P_{(1,1,-1,-1)}= & \Psi_{(3,3)/(2,2)}(1-t^{\kappa}q)
  (1-t^{\kappa-1}q) P_{(2,2,-1,-1)}\\  +&
  \Psi_{(3,2,1)/(2,2)}(1-t^{\kappa}q) (1-t^{\kappa-2}q^{-1})
  P_{(2,1,0,-1)}\\ 
  + & \Psi_{(2,2,1,1)/(2,2)}(1-t^{\kappa-2}q^{-1})
  (1-t^{\kappa-3}q^{-1}) P_{(1,1,0,0)}
\end{split}
\end{equation}

Remarkably it seems that the creation operators $\tilde B_k^+$ can be
identified with a subset of the operators $F_{m,\kappa}$.  Indeed the
following conjecture has been arrived at with the aid of the computer.
\begin{conjecture} The creation operators $\tilde B_k^+$ can be
  written in the form 
\begin{equation}
\tilde B_k^+= F(k) e_k F(k)^{-1}=F_{k,k}
\end{equation}
\end{conjecture}
This expression immediately provides, through (60), the action of
$\tilde B_k^+$ on arbitrary Macdonald polynomials.  Conjecture~4 must
then be a consequence of it.  To convince oneself that formula (64)
indeed implies Conjecture~4, one uses the same kind of argument as in
the proof of Lemma~2. Evaluating the action of $F_{k,k}$ on
$P_{\lambda}$ with the help of (60), one thus shows that all the terms
associated to partitions with more than $k$ parts are annihilated.  In
the framework of this conjecture, the Hermitian conjugate $\tilde
B_k^-$ of $\tilde B_k^+$ with respect to the scalar product defined in
(22) are readily obtained from the fact that $F(\kappa)$ is Hermitian
under this scalar product.  We thus have
\begin{equation}
F_{m,\kappa}^{\dagger}= F(\kappa)^{-1}  e_m^{\dagger} F(\kappa)
\end{equation}
with $ e_m^{\dagger}= e_N^{-1} e_{N-m}$, which implies that
\begin{equation}
\tilde B_k^{-}= (\tilde B_k^+)^{\dagger} =  F(k)^{-1}  e_k^{\dagger} F(k).
\end{equation}
It is striking that the set of operators $F_{m,\kappa,\gamma}$
contains a one-parameter family of $N$-dimensional Abelian algebras.
From (58), all $F_{m,\kappa}$ can be generated from $F_{m,0}$ by
conjugating with powers of $e_N$.  More precisely, one has
\begin{equation}
F_{m,\kappa}= e_N^{-\kappa/\alpha}  F_{m,0} e_N^{\kappa/\alpha} .
\end{equation}
With
\begin{equation}
[ F_{m,\kappa},F_{n,\kappa}]=F(\kappa)[e_m,e_n] F(\kappa)^{-1}=0,
\end{equation}
this immediately shows that we can construct, by proper conjugation of
$\tilde B_k^+$ with $e_N$, the following set of commuting operators:
\begin{equation}
\{ F_{m,\kappa}= e_N^{(m-\kappa)/\alpha} \tilde B_m^+
e_N^{(\kappa-m)/\alpha} , m=1,\dots,N ;\kappa \in \mathbb R \}  
\end{equation}
and, for each value of $\kappa$, thus obtain a $N$-dimensional Abelian
algebra.

A straightforward extension of Lemma~6 is
\begin{equation}
M_I(X;q,t) e_N^{\rho} =  e_N^{\rho} M_I(X q^{\rho};q,t),
\end{equation}
which can in particular be seen from the fact that $M_I^r
e_N^{\rho}=(q^r)^{\rho} e_N^{\rho} M_I^r$.  This last result finally
provides the following realization of the operators $ F_{m,\kappa}$:
\begin{equation}
 F_{m,\kappa}= \sum_I \tilde A_I(x;t) x_I M_I(-t^{\kappa-m+1};q,t)
\end{equation}
assuming that Conjecture~8 is valid.

\section{Jack polynomials \cite{1}}
The monic Jack polynomials $P_{\lambda}(x;\alpha)$ are obtained from
the monic Macdonald polynomials $P_{\lambda}(x;q,t)$ in the limit
\begin{equation}
q=t^{\alpha},\quad t \to 1.
\end{equation}
Let $\mathbb Q(\alpha)$ denote the field of rational functions of
$\alpha$.  Taking the above limit in (7) yields the following scalar
product on $\Lambda_N \otimes \mathbb Q(\alpha)$:
\begin{equation}
\langle p_\lambda, p_\mu \rangle_{\alpha}
        =\delta_{\lambda \mu } z_\lambda \alpha^{\ell(\lambda)},
\end{equation}
with $\ell(\lambda)$, the number of parts of $\lambda$.  The Jack
polynomials $J_{\lambda}(x;\alpha)$ $\in \Lambda_N \otimes \mathbb
Q(\alpha)$ are given by
\begin{equation}
J_\lambda (x;\alpha)=\lim_{\begin{subarray}{c} q=t^{\alpha} \\ t \to 1
\end{subarray}} \frac{J_{\lambda}(x;q,t)}{(1-t)^{|\lambda|}},
\end{equation}
with $|\lambda|=\lambda_1+\lambda_2+\dots$. They are thus uniquely
specified by
\begin{align}
  \mathrm{(i)} \ &  \langle J_\lambda, J_\mu \rangle_{\alpha} = 0, \qquad
  \text{if } \lambda \ne \mu , \label{3}\\ 
  \mathrm{(ii)} \ &  J_\lambda = \sum_{\mu \le \lambda}
  v_{\lambda\mu }(\alpha) m_\mu ,  \label{4}\\ 
  \mathrm{(iii)} \ & v_{\lambda\lambda}(\alpha)= c_{\lambda}(\alpha),
\label{5}
\end{align}
where
\begin{equation}
c_{\lambda}(\alpha) = \lim_{\begin{subarray}{c} q=t^{\alpha} \\ t \to 1
\end{subarray}} \frac{c_{\lambda}(q,t)}{(1-t)^{|\lambda|}}=\prod_{s
\in \lambda} (\alpha a(s) + \ell (s) +1). 
\end{equation}
In the limit (78) the Pieri formula becomes
\begin{equation}
e_k P_{\lambda}^{(\alpha)} = \sum_\mu  {\Psi}_{\mu /\lambda}^{(\alpha)}
P_\mu ^{(\alpha)}, 
\end{equation}
with the sum over partition $\mu \supset \lambda$ (of length $\le N$)
such that $\mu -\lambda$ is a vertical $k$-strip.  The coefficients
${\Psi}_{\mu /\lambda}^{(\alpha)}$ are
\begin{equation} 
{\Psi}_{\mu /\lambda}^{(\alpha)} = \prod_{\begin{subarray}{c} s \in
    C_{\mu /\lambda} \\ s  
\not \in R_{\mu /\lambda} \end{subarray}}
\frac{b_\mu ^{(\alpha)}(s)}{b_{\lambda}^{(\alpha)}(s)}, 
\end{equation}
where
\begin{equation}
b_{\lambda}^{(\alpha)}(s) = {\cases \frac{\alpha a(s)
    +\ell(s)+1}{\alpha(a(s)+1)+\ell(s)} \qquad & {\text{if}}~s \in
  \lambda \\ 
1    & {\text{otherwise}} \endcases}.
\end{equation}

We shall now recall the definition of the creation operators entering
in the Rodrigues formula for the Jack polynomials that we derived in
\cite{3} and shall establish their connection with the operators
$\tilde B_k^+$ in the limit (72).  The creation operators in this case
are constructed from the Dunkl operators
\begin{equation} 
D_i
        = \alpha x_i \frac{\partial}{\partial x_i} +
        \sum_{\begin{subarray}{l}j=1\\ j\ne i \end{subarray}}^N
        \frac{x_i}{x_i - x_j} (1 - K_{ij}), \qquad i = 1,2, \dots , N,
\end{equation}
where $K_{ij} = K_{ji}$ is the operator that permutes the variables
$x_i$ and $x_j$:
\begin{equation} 
K_{ij}x_i
        = x_j K_{ij}, \quad K_{ij} D_i = D_j K_{ij}, \quad K_{ij}^2 = 1.
\end{equation}
Let $J = \{ j_1, j_2, \dots , j_\ell \}$ be sets of cardinality $|J| =
\ell$ made of integers $j_\kappa \in \{ 1, \dots , N \}$, $1 \le
\kappa \le \ell$ such that $j_1 < j_2 < \dots < j_\ell$ and introduce
the operators
\begin{equation} 
D_{J,\omega}
        = (D_{j_1} + \omega) (D_{j_2} + \omega + 1) \dotsm (D_{j_\ell}
        + \omega +\ell-1), 
\end{equation}
labelled by such sets and a real number $\omega$.  The creation
operators $\tilde B_k^{+ (\alpha)}$ are defined by
\begin{equation} 
\tilde B_k^{+ (\alpha) }
        = \sum_{\begin{subarray}{c} J \subset \{ 1,\dots , N\}\\ |J| =
            k \end{subarray}}  x_J D_{J,1}. 
\end{equation}
From the following theorem proved in \cite{2}, we see that the
operators $\tilde B_k^{+ (\alpha) }$ construct the Jack polynomials.
\begin{theorem} The Jack polynomials $J_\lambda(x;\alpha)$ associated
  to the partitions $\lambda$ $=$ $(\lambda_1$ , $\dots$ ,
  $\lambda_N)$ are given by 
\begin{equation} 
J_\lambda(x;\alpha)
        = (\tilde B_N^{+(\alpha)})^{\lambda_N} (\tilde
        B_{N-1}^{+(\alpha)})^{\lambda_{N-1} -\lambda_N} \dots (\tilde
        B_1^{+(\alpha)})^{\lambda_1 - \lambda_2} \cdot 1 \quad . 
\end{equation}
\end{theorem}
We shall now show that the operators $\tilde B_k^{+(\alpha)}$ entering
in this formula are related to the $\tilde B_k^+$ of Conjecture~4 and
hence that the conjectures of section~4 on the $\tilde B_k^+$ also
apply for the operators $\tilde B_1^{+(\alpha)}$.
\begin{lemma}  The $\tilde B_k^{+ (\alpha) }$ are the following limits
  of the $\tilde B_k^{+ }$ 
\begin{equation}
\Res \tilde B_k^{+(\alpha)}=\lim_{\begin{subarray}{c} q=t^{\alpha} \\ t \to 1
\end{subarray}} \frac{\tilde B_k^+}{(1-t)^k} \quad .
\end{equation}
\end{lemma}
Here $\Res X$ means that $X$ is restricted to symmetric functions of
the variables $x_1,\dots,x_N$.  It has been shown in \cite{2} that
$\Res D_{J,\omega}$ is symmetric under the permutations of the
variables $x_i, i\in J$ and depends only upon these variables.  By an
argument such as the one given in section 4.~4 of \cite{2}, we have
that
\begin{equation}
\Bigl( \Res D_{J,\omega} \Bigr) J_{\lambda}(x(J);\alpha)=
\prod_{i=1}^{\ell} (\alpha \lambda_i +\omega+\ell-i)
J_{\lambda}(x(J);\alpha), 
\end{equation}
with $x(J)=\{ x_i | i \in J \}$.  Moreover,
\begin{equation}
\begin{split}
  \lim_{\begin{subarray}{c} q= t^{\alpha} \\ t \to 1
\end{subarray}} \Bigl( \frac{M_J(-t^{\omega};q,t)}{(1-t)^{\ell}}
\frac{J_{\lambda}(x(J);q,t)}{(1-t)^{|\lambda|}} \Bigr)  
=& \Biggl( \lim_{\begin{subarray}{c} q= t^{\alpha} \\ t \to 1
\end{subarray}} \prod_{i=1}^{\ell} \frac{(1-q^{\lambda_i} t^{\omega+\ell-i})}
{(1-t)} \Biggr) J_{\lambda}(x(J);\alpha)\\
=& \prod_{i=1}^{\ell} (\alpha \lambda_i +\omega+\ell-i)
J_{\lambda}(x(J);\alpha)\\ 
=& \Bigl( \Res D_{J,\omega} \Bigr) J_{\lambda}(x(J);\alpha).
\end{split}
\end{equation}
Since $J_{\lambda}(x(J);\alpha)$ form a basis for the symmetric
polynomials in the variables $x_i, i \in J$, we have that
\begin{equation}
\lim_{\begin{subarray}{c} q=t^{\alpha} \\ t \to 1
\end{subarray}} \frac{M_J(-t^{\omega};q,t)}{(1-t)^{\ell}}= \Res D_{J,\omega},
\end{equation}
knowing that differential operators that coincide on a complete set of
symmetric functions are identical (see for instance appendix C of
\cite{5}).  We thus have
\begin{equation}
\lim_{\begin{subarray}{c} q=t^{\alpha} \\ t \to 1
\end{subarray}} \frac{\tilde B_k^+}{(1-t)^k}= \sum_{|J|=k} x_J \Res D_{J,1} = 
\Res \sum_{|J|=k} x_J D_{J,1},
\end{equation}
which proves Lemma~10.

In light of Theorem~1 and (74) , we can write in addition to (86)
another Rodrigues formula for the Jack polynomials using
$\lim_{\begin{subarray}{c} q=t^{\alpha} \\ t \to 1
\end{subarray}} B_k^+/(1-t)^k $ as creation operators.  With the help
of  (90), one gets 
\begin{equation}
\begin{split}
  B_k^{+(\alpha)} = \lim_{\begin{subarray}{c} q=t^{\alpha} \\ t \to 1
\end{subarray}} B_k^+/(1-t)^k = &\lim_{\begin{subarray}{c}
  q=t^{\alpha} \\ t \to 1 \end{subarray}} \prod_{j=k+1}^N \Biggl[
\frac{(1-t)}{(1-t^{k+1-j}q^{-1})} \Biggr]
\frac{M_N(-t^{k+1-N}q^{-1};q,t)}{(1-t)^N} e_k \\ 
=& \prod_{j=k+1}^N (-\alpha +k+1-j)^{-1}
D_{\{1,\dots,N\},k+1-N-\alpha}~e_k.
\end{split}
\end{equation}
These operators $B_k^{+(\alpha)}$ are also such that
\begin{equation} 
J_\lambda(x;\alpha)
        = ( B_N^{+(\alpha)})^{\lambda_N} (
        B_{N-1}^{+(\alpha)})^{\lambda_{N-1} -\lambda_N} \dots (
        B_1^{+(\alpha)})^{\lambda_1 - \lambda_2} \cdot 1,  
\end{equation}
for any partition $\lambda$.

As in the case of the Macdonald polynomials, let us extend the
definition of a partition to allow real entries:
\begin{equation}
P_{(\beta_1,\dots,\beta_N)}^{(\alpha)} = e_N^{\beta_N}
P_{(\beta_1-\beta_N,\dots,\beta_{N-1}-\beta_N,0)}^{(\alpha)}=
e_N^{\beta_N} P_{\beta-\beta_N}^{(\alpha)} 
\end{equation}
$\forall \beta_N \in \mathbb R$ and $\beta_{i}-\beta_{i+1}$ a
non-negative integer, $i=1,\dots,N-1$.  Take the operator
$F^{(\alpha)}(\kappa)$ that acts as follows on $P_{\beta}^{(\alpha)}$:
\begin{equation}
F^{(\alpha)}(\kappa)
P_{(\beta_1,\dots,\beta_N)}^{(\alpha)}=\prod_{i=1}^N
\prod_{j=1}^{\infty} (\alpha(\beta_i-j)+\kappa+1-i)
P_{(\beta_1,\dots,\beta_N)}^{(\alpha)}, 
\end{equation}
Again, it follows that
\begin{equation}
F^{(\alpha)}(\kappa) e_N^{\rho}=e_N^{\rho} F^{(\alpha)}(\kappa+\alpha \rho).
\end{equation}
We now form the operators
\begin{equation}
F_{m,\kappa}^{(\alpha)} = F^{(\alpha)}(\kappa) e_m F^{(\alpha)}(\kappa)^{-1},
\end{equation}
to see that these have on $P_{\beta}^{(\alpha)}$ actions that only
involve a finite number of products.  These actions read
\begin{equation}
F_{m,\kappa}^{(\alpha)} P_{\beta}^{(\alpha)} = \sum_{\delta}
{\Psi}_{\delta/\beta}^{(\alpha)}  
F_{\delta/\beta}^{(\alpha)}(\kappa) P_{\delta}^{(\alpha)},
\end{equation}
with $\delta-\beta$ $m$-vertical strips,
\begin{equation}
F_{\delta/\beta}^{(\alpha)}(\kappa)= \prod_{\begin{subarray}{c} s \in
    \delta-\beta_N \\ s \not \in \beta-\beta_N 
\end{subarray}} F_{\delta-\beta_N}^{(\alpha)}(s;\kappa+\alpha \beta_N),
\end{equation}
and where for partitions $\mu $ made of nonnegative integers,
\begin{equation}
F_\mu ^{(\alpha)}(s;\kappa)=\bigl( \alpha(a'(s))+\kappa-\ell'(s)
\bigr),\quad \forall s \in \mu . 
\end{equation}
Again, ${\Psi}_{\delta/\beta}^{(\alpha)} =
{\Psi}_{\delta-\beta_N/\beta-\beta_N}^{(\alpha)}$.  That
\begin{equation}
F_{m,\kappa}^{(\alpha)}= \lim_{\begin{subarray}{c} q= t^{\alpha} \\ t
    \to 1 \end{subarray}} \frac{F_{m,\kappa}}{(1-t)^m}, 
\end{equation}
is established from the fact that
\begin{equation}
\lim_{\begin{subarray}{c} q= t^{\alpha} \\ t \to 1 \end{subarray}}
\frac{F_{\delta/\beta}(\kappa)}{(1-t)^m}=F_{\delta/\beta}^{(\alpha)}(\kappa),
\end{equation}
with $\delta-\beta$ a vertical m-strip. It then follows that
Conjecture~8 implies, if true, that the operators $\tilde
B_k^{+(\alpha)}$ and $\tilde B_k^{+}$ share many features.
\begin{conjecture} The creation operators $ \tilde B_k^{+(\alpha)}$
  share these properties: 
\begin{equation}
\begin{split}
  \mathrm{(i)} \ & \tilde  B_k^{+(\alpha)}= F^{(\alpha)}(k) e_k
  F^{(\alpha)}(k)^{-1}=F_{k,k}^{(\alpha)}\\ 
  \mathrm{(ii)} \ &  F_{m,\kappa}^{(\alpha)}= \sum_{|J|=m} x_J
  D_{J,\kappa-m+1} = e_N^{(m-\kappa)/\alpha } \tilde B_m^{+(\alpha)}
  e_N^{(\kappa-m)/\alpha }\\ 
  \mathrm{(iii)} \ & [ F_{m,\kappa}^{(\alpha)},
  F_{n,\kappa}^{(\alpha)}]=0 \quad.
\end{split}
\end{equation}
\end{conjecture}
In view of (98), we thus have a conjecture for the action of $\tilde
B_k^{+(\alpha)}$ on arbitrary Jack polynomials.

\begin{acknow}
  We would like to express our thanks to Fran{\c{c}}ois Bergeron and
  Adriano Garsia for various comments and suggestions.  Extensive use
  was made of J. Stembridge's Symmetric Functions Maple package.

\noindent This work has been supported in part through funds provided
by NSERC (Canada) and FCAR (Qu{\'e}bec).  L.~Lapointe holds a NSERC
postgraduate scholarship. 
\end{acknow}

     \end{document}